\documentclass[11pt]{amsart}
\usepackage{amsmath,amssymb, amsthm}

\newcommand\N{{\mathbb N}}

\newcommand\G{\Gamma}
\newcommand\inv{^{-1}}

\newtheorem{theorem}{Theorem}[section]

\newtheorem{lemma}[theorem]{Lemma}
\newtheorem{corollary}[theorem]{Corollary}

\newtheorem{example}[theorem]{Example}
\newtheorem{question}[theorem]{Question}

\title{A note on spaces of asymptotic dimension one}

\author{Koji Fujiwara}
\address{Math Institute, Tohoku University, Sendai,
980-8578 Japan}
\email{fujiwara@math.tohoku.ac.jp}
\thanks{This work has been done when the first author was 
visiting the second author at University of Illinois
at Chicago, and during his visit 
at Max-Planck-Institute for Mathematics in Bonn.
He would like to thank both institutions for hospitality.}

\author{Kevin Whyte}
\address{Department of Mathematics, University of Illinois at Chicago, Chicago, Il }
\email{kwhyte@math.uic.edu}
\thanks{The second author was partially supported by an NSF Career award (DMS-0204576) }


\date{2006.11.29}

\begin{document}

\renewcommand{\thefootnote}{\fnsymbol{footnote}}
\footnote[0]{2000\textit{ Mathematics Subjet Classification}.
51F99; 20F69, 54F45, 57M50}

\maketitle

\begin{abstract}
Let $X$ be a geodesic metric space with $H_1(X)$ uniformly generated.  
If $X$ has asymptotic dimension one then $X$ is quasi-isometric to an unbounded tree.
As a corollary, we show that the asymptotic dimension of
 the curve graph of a compact, oriented surface 
with genus $g \ge 2$ and one boundary component
 is at least two.
\end{abstract}

\section*{Introduction}

Our main results show that, under fairly weak restrictions, a geodesic metric space of asymptotic dimension one is quasi-isometric to a tree.   Before stating precise results we review some definitions.

If $X$ is a set and $X=\cup_i O_i$ a covering, we say that the {\it multiplicity} of the covering 
is at most $n$ if any point $x \in X$ is contained in at most $n$ elements of $\{ O_i\}_i$.  

Recall that the {\it covering dimension} of a topological space $X$ is the minimal $d$ such that  every open covering has a refinement with multiplicity at most $d+1$.   The notion of asymptotic dimension is a coarse analogue introduced by Gromov in \cite{Gr93}: 

\noindent
Let $X$ be a metric space, and $X=\cup_i O_i$
a covering. For $D \ge 0$, we say 
that the {\it $D$-multiplicity} of the covering 
is at most $n$ if for any $x \in X$, 
the closed $D$-ball centered 
at $x$ intersects at most $n$ elements of
$\{ O_i\}_i$. The multiplicity is exactly the
$0$-multiplicity.

The {\it asymptotic dimension} of the metric space $X$ 
is at most $n$ if for any $D\ge 0$, there
exists a covering $X=\cup_i O_i$ such 
that the diameter of $O_i$ is
uniformly bounded (i.e. there exists $C$ such that 
for all $i$, ${\rm diam} O_i \le C$), and 
the $D$-multiplicity of the covering is at most $n+1$.
We say that the asymptotic dimension of $X$, $asdim \, X$, is $n$
if the asymptotic dimension of $X$ is at most $n$,
but it is not at most $n-1$. If such $n$ does not exist, 
then we define the asymptotic dimension of $X$ to be infinite.

It is easy to see that if two metric spaces
are quasi-isometric, then they have the same asymptotic dimension.
A geodesic metric space has asymptotic dimension zero
if and only if it is bounded.
Since unbounded trees have asymptotic dimension one, all spaces quasi-isometric to 
unbounded trees also have asymptotic dimension one, \cite{Gr93}.   

Let $X$ be a geodesic metric space.  We say $H_1(X)$ is {\it uniformly generated} if there is an $L>0$ so that $H_1(X)$ is generated by loops of length at most $L$.  

\begin{theorem}\label{main}  Let $X$ be a geodesic metric space with $H_1(X)$ uniformly generated.  
If $X$ has asymptotic dimension one then $X$ is quasi-isometric to an unbounded tree.
\end{theorem}

\noindent
{\bf Remark}: This theorem relies on a result of Manning (\cite{Man}) characterizing spaces quasi-isometric to trees, see Theorem \ref{qtrees}.  The restriction on $H_1(X)$ is essential, see example \ref{pearls}.

Let $G$ be a finitely generated group, and $\Gamma$
the Cayley graph with respect to a finite generating set.
Since the asymptotic dimension is invariant by quasi-isometry,
the asymptotic dimension of $G$, $asdim \, G$, is defined to be the asymptotic 
dimension of $\Gamma$.

\begin{corollary}\label{fp}
Let $G$ be a finitely presented group. The asymptotic dimension of $G$
is one if and only if  $G$ contains a non-trivial free group as a subgroup of finite index.
\end{corollary}

We are informed by Dranishnikov that this has been known \cite{JS} (cf. \cite{Ge}, \cite{D}).
The assumption that $G$ is finitely presented is necessary.
The following example is due to Osin.
\begin{example}
Let $G=A \wr {\Bbb Z}$ be the wreath product such that $A$ is a non-trivial 
finite group. Then, ${\rm asdim} G=1$ but $G$ does not contain a free group 
as a subgroup of finite index. $G$ is finitely generated, but not finitely 
presented.
\end{example}
By definition, $G$ is the semidirect product of $A^{{\Bbb Z}}$ and ${\Bbb Z}$,
$A^{{\Bbb Z}} \rtimes {\Bbb Z}$,
such that the action of ${\Bbb Z}$ is the obvious action of shifting 
the indexes of the direct product, $N=A^{{\Bbb Z}}$, of countable copies of $A$.
There is an exact sequence $1 \to N \to G \to {\Bbb Z} \to 0$.
Note that $N$ is a locally finite (i.e. any finitely generated subgroup is finite)
countable group, which is not finitely generated.
Therefore $G$ does not contain a free group as a subgroup of finite index.
We want to show that 
${\rm asdim} G=1$. Indeed, the notion of asymptotic dimension 
is extended to a countable group, and it is shown that 
${\rm asdim} G= \sup {\rm asdim} F$, where sup is taken over all finitely 
generated subgroups $F<G$, \cite{DS}.
Since $N$ is locally finite, 
${\rm asdim} N=0$. The Hurewicz type formula for asymptotic dimension 
is also extended, \cite{DS}, 
and we get 
${\rm asdim}G \le {\rm asdim} N+{\rm asdim}{\Bbb Z}=1$.
Also, since ${\Bbb Z} < G$, $1 \le {\rm asdim}G$. 
It is easy to see that if $K$ and $L$ are finitely generated, then 
the wreath product $K \wr L$ is finitely generated.
If both $K,L$ are finitely presented, 
$K \wr L$ is finitely presented if and only if 
$K$ is trivial or $L$ is finite, \cite{Ba}. 
Therefore, $G$ is finitely generated, but 
not finitely presented.

It is natural to ask the following question:
\begin{question}
Suppose $G$ is a finitely generated group which 
is torsion free. If ${\rm asdim} G=1$, then is $G$ a free group ?
\end{question}

\begin{corollary}\label{hyp}
Let $X$ be geodesic metric space which is Gromov hyperbolic.  
If $X$ has asymptotic dimension one then $X$ is quasi-isometric to an unbounded tree.
\end{corollary}

This result does not require that $X$ be proper.  One important class of examples which are not proper are  curve complexes. Let $S=S_{g,p}$ be a compact, orientable surface
such that $g$ is the genus and $p$ is the 
number of the connected components of the boundary of $S$.  We assume
that $3g-4+p >0$.
The {\it curve complex} of $S$, defined by
Harvey \cite{Ha}, are the flag complexes with 1-skeleton the {\it curve graph} of $S$, $C(S)$.   
The curve graph is the graph whose vertices are isotopy classes of
essential, nonperipheral, simple closed
curves in $S$, with two distinct
vertices  joined by an edge if the corresponding curves can
be realized by disjoint curves. 
We remark that the curve complex of $S$ is quasi-isometric
to the curve graph of $S$, so that they have same 
asymptotic dimension.

Masur and Minsky \cite{MaMi} show the remarkable result
that $C(S)$ is Gromov hyperbolic.

It is known that $C(S)$ has finite asymptotic dimension (\cite{BeF}).  
Since $C(S)$ is unbounded  it cannot have asymptotic dimension zero.  
No upper bound or non-trivial lower bound is known in general.   
In the case $p=1$, Schleimer shows in \cite{Sch} that $C(S)$ is one-ended, 
from which it follows that it is not quasi-isometric to a tree.  
Thus we can improve the lower bound on asymptotic dimension here:

\begin{corollary}
Let $S$ be a compact, oriented surface 
with genus $g \ge 2$ and one boundary component.
Let $C(S)$ be the curve graph of $S$.
Then the asymptotic dimension of $C(S)$ is 
at least two.
\end{corollary}

\noindent
{\bf Acknowledgements.}
The first author would like to thank Greg Bell for 
discussion. We appreciate comments and information by Alexander Dranishnikov.
We are grateful to Denis Osin, who shared his ideas with us.

\section{Proofs}

Theorem \ref{main} follows from two key ingredients.   The first is the classical argument that surfaces have covering dimension two (see, for instance, \cite{Mu}, section 55).  The fact we need is the following:

\begin{lemma}\label{cover} Let $S$ be a compact surface with boundary $S^1$.   
Let $A$, $B$, and $C$ be points on $\partial S$, dividing it into arcs $AB$, $AC$, and $BC$.   
 If $\{ O_i\}_i$ is an open cover of $S$ with multiplicity two then there is some $i$ for which $O_i$ intersects all three segments $AB$, $AC$, and $BC$.
\end{lemma}

The second ingredient is a characterization of spaces quasi-isometric to trees due to Manning:

\begin{theorem}[Th 4.6 \cite{Man}]\label{qtrees}
Let $Y$ be a geodesic space.
Suppose that there exists a constant $K$ with 
the following property:\\
\\
let $a,b$ be points in $Y$, $\gamma$
a geodesic from $a$ to $b$, and $\alpha$ a path 
from $a$ to $b$. Then the $K$-ball at the midpoint 
$m$ of the geodesic $\gamma$ has non-empty intersection 
with $\alpha$. \\

Then $Y$ is quasi-isometric to a simplicial tree.
\end{theorem}

\noindent
We now prove Theorem \ref{main}:

\begin{proof}

The following lemma, essentially contained in \cite{BW},  explains the role of the uniform generation assumption:

\begin{lemma}  Let $X$ be a complete geodesic metric space.   The following are equivalent:
\begin{itemize}
\item $X$ has uniformly generated $H_1$.
\item $X$ is quasi-isometric to a complete geodesic metric space,  $Y$, with $H_1(Y)=0$.
\end{itemize}
\end{lemma}

\noindent
(These are also equivalent to the condition $H_1^{uf}(X)=0$ which is the form which appears in \cite{BW}.)

\begin{proof}

Let $Y$ be a complete geodesic metric space with $H_1(Y)=0$, and let $f: X \to Y$ be a $(K,C)$-quasi-isometry.      Suppose $l$ is any loop in $X$.     Divide $l$ into segments $[a_0,a_1],[a_1,a_2],\ldots , [a_{n-1},a_n]$, where $a_0=a_n$ and the length of each segment is at most $1$.    

Let  $b_i = f(a_i)$.  We have $d(b_i,b_i+1) \leq K+C$.  For each $i$ choose a geodesic segment in $Y$ connecting $b_i$ to $b_{i+1}$.   Call the resulting loop $l'$.    By assumption, $l'$ bounds, so we have $l'= \partial S$ for some surface $S$ mapping to $Y$.   Triangulate $S$ so that the image of each simplex has diameter at most $1$ in $Y$.     For each vertex $v$ in $S$ let $\alpha(v)$ be a point of $X$ which maps close to the image of $v$ in $Y$.    For each edge $e$ of $S$ with endpoints $u$ and $v$, choose a geodesic segment, $\alpha(e)$, in $X$ connecting $\alpha(u)$ to $\alpha(v)$. 

Thus we have $\alpha$ mapping the one-skeleton of $S$ to $X$.  For each simplex $\sigma$ of $S$ we have a loop in $X$ of length at most $D$ (depending only on $K$ and $C$).     
The boundary component of $S$ is mapped to 
a loop $l''$ which is within Hausdorff distance $D$ of $l$.    The map $\alpha$ exhibits that $l''$ is a sum of classes of length at most $D$, and the condition on Hausdorff distance implies the same for $l - l''$.   Thus $l$ is a sum of loops of length at most $D$, proving that $H_1(X)$ is uniformly generated.  

Conversely, suppose $H_1(X)$ is generated by loops of length at most $D$.    Let $A$ be a maximal collection of points in $X$ with $d(a,a') \geq D$ for all $a\neq a'$.     Let $R=3D$, and let $Y$ be the space: 
$$ X  \cup_{a \in A} cone(B(a,R)) $$
\noindent

In words, $Y$ is $X$ with each $R$-ball centered at $a \in A$ coned to a point.   We give $Y$ the induced path metric where each cone line has length $R$.   The inclusion of $X$ is then isometric, and has coarsely dense image, and so is a quasi-isometry.  By construction, any loop $l$ in $X$ of length at most $D$ is contained $B(a,R)$ for some $a$, and hence is null homotopic in $Y$.  Since these generate $H_1(X)$, we have $H_1(X) \to H_1(Y)$ is the zero map.     

Let $l$ be any loop in $Y$.   By compactness of the circle,  $l$ only passes through finitely many cone points.    By pushing $l$ off of the cones whose cone points it misses, we may arrange for $l$ to visit the interior of only finitely many cones.   For each such cone $cone(B(a,R))$, let $u$ and $v$ be the first and last places $l$ visits this cone.    By a further homotopy we may assume $l$ travels along the cone line from $u$ to the cone point and then along the cone line back from the cone point to $v$.     Choose a path $p_a$ of length at most $2R$ in $B(a,R)$ connecting $u$ and $v$.      Let $l'$ be the loop with all the trips into the cones cut out and replaced with the $p_a$.  The disks $c(p_a)$ give a homotopy between $l$ and $l'$ in $Y$. Since $l' \subset X$ we know $[l]=[l']=0$ in $Y$.    Thus $H_1(Y)=0$ as claimed.
\end{proof}

\noindent
{\bf Remark}: The above proof goes through verbatim to show that uniform generation of $\pi_1(X)$ is equivalent to $X$ quasi-isometric to a simply connected space.

\noindent
Returning to the proof of Theorem \ref{main}, we assume $X$ is
quasi-isometric to a geodesic space $Y$ with  $H_1(Y)=0$ using the lemma.
Since $Y$ has asymptotic dimension one there is a cover $O_i$ of $Y$ with $1$-multiplicity two, 
and all $O_i$ of diameter at most some $D < \infty$.   
Let $U_i$ be the open subset of $Y$ consisting of the points within a distance of $1$ of $U_i$.  
By definition this cover has multiplicity at most two, and  each $U_i$ has diameter at most $D + 2$.

We check that Manning's condition holds in $Y$ with $K=\frac{3}{2}(D+2)$.   
Let $\gamma$, $a$, $b$, $m$, and $\alpha$ be as in the statement of theorem \ref{qtrees}.     
Let $l$ be the loop $\gamma \circ \alpha\inv$.    
Since $H_1(Y)=0$ there is a surface $S$ with $\partial S = S^1$ and a continuous map $f:S \to Y$ 
with $f|_{\partial S} = l$.  Applying lemma \ref{cover} to the covering $f\inv (U_i)$ shows that 
there is some $i$ such that $U_i$ has non-trivial intersection with $\alpha$ and with 
the segments $\overline{am}$ and $\overline{mb}$ of $\gamma$.    
Since $U_i$ has diameter at most $D+2$ this means there are points $p$ on $\alpha$ 
and $s$ and $t$ on $\overline{am}$ and $\overline{mb}$ which are pairwise within a distance $D+2$.    
Since $\gamma$ is a geodesic, one of $s$ or $t$ must be within $\frac{1}{2}(D+2)$ of $m$, 
so $\alpha$ must pass through the $K$ ball around $m$ as desired.
We showed that $Y$, therefore $X$ as well, is quasi-isometric to a tree, which 
is unbounded.
\end{proof}

\noindent
{\bf Remark}: The proof only uses that $X$ has a cover by sets of bounded diameter with $L_0$-multiplicity two for $L_0$ large enough to generate all of $H_1(X)$.  Since it follows that $X$ is quasi-isometric to a tree, 
we know that there are covers by uniformly bounded sets of $L$-multiplicity two for all $L$.\\

\noindent
Theorem \ref{main} does not hold without the assumption on $H_1(X)$:\\

\noindent
\begin{example}\label{pearls}
{
Let $X$ be a graph with vertex set $\N$ and two edges connecting $n$ and $n+1$, both of length $2^n$.   Thus $X$ is a string of circles of increasing radii, each touching the next at a single point.   It is clear that $H_1(X)$ is not uniformly generated and in particular that $X$ is not quasi-isometric to a tree.  

For any $D>0$, choose $n$ such that $2^n >> D$.   Let $U$ be the union of the vertices $\{1,2, \cdots, n\}$ of $X$ and all adjacent edges except for the $\frac{1}{2}D$-ball around the vertex $n+1$.   For $i>0$ let $A_i$ be the ball of radius $D$ around the vertex $n+i$, let $B_i$ and $C_i$ be the subsets of the two edges from $n+i$ to $n+i+1$ not within $\frac{1}{2}D$ of a vertex.    It is easy to see that collection $\{U,A_i,B_i,C_i\}$ is a covering of $X$ with $D$-multiplicity two.   By subdividing the $B_i$ and $C_i$ to a collection of intervals of length $2D$ covering each, one gets a covering with $D$-multiplicity by sets of diameter at most $2^n$.   Thus $X$ has asymptotic dimension one.
}
\end{example}

\noindent
{\it Proof of Cor \ref{fp}}.
A well known corollary of work of Stallings (\cite{Sta}) and Dunwoody (\cite{Dun}) is that 
a group quasi-isometric to a tree is virtually free.   
Thus corollary \ref{fp} follows immediately from Theorem \ref{main} 
since the fundamental group of the Cayley graph is generated by translates of the relators, 
and so is generated by loops of length at most the length of the longest relator.   
\qed

\vspace{3mm}
\noindent
{\bf Remark}: 
 In fact, we get something a priori stronger: if $\G = \langle F | R \rangle $ where $R_{ab}$ 
is a finitely generated $\G$-module and $\G$ has asymptotic dimension one then $\G$ 
is virtually free.   However, since every virtually free group is finitely presented, 
this does not include any new examples.  

\vspace{3mm}
\noindent
{\it Proof of Cor \ref{hyp}}.
To prove Corollary \ref{hyp}, we need to check that $H_1(X)$ is uniformly generated.  
For any graph $G$ we have that $H_1(G)$, and, indeed, $\pi_1(G)$, are generated by 
isometrically embedded loops.   
If we divide such a loop into three equal pieces we get a geodesic triangle 
whose sides come close only at the corners.   
If $G$ is $\delta$-hyperbolic, this implies any such loop can have length at most $6\delta$, 
which proves the corollary.   
\qed

\medskip
\noindent


\begin{thebibliography}{Dun85}

\bibitem[Ba]{Ba}
G.Baumslag.
Wreath product and finitely presented groups.
Math Zeit. 75, (1961). 22-28. 

\bibitem[BeF]{BeF}
G.Bell, K.Fujiwara.
The asymptotic dimension of a curve graph is finite,
preprint, 2005 September. Arxiv, math.GT/0509216.

\bibitem[BW]{BW}
J. Block and S. Weinberger, {\em Large Scale Homology Theories and 
Geometry}, AMS/IP Studies in Advanced Mathematics, volume 2, 1997, pp. 522-569.

\bibitem[D]{D}
A. N. Dranishnikov.
Cohomological approach to asymptotic dimension. 
preprint, 2006 August, math.MG/0608215.
 
\bibitem[DS]{DS}
A.Dranishnikov, J.Smith.
Asymptotic dimension of discrete groups. 
Fund. Math. 189 (2006), no. 1, 27--34. 

\bibitem[Dun85]{Dun}
M.~J. Dunwoody, \emph{The accessibility of finitely presented groups}, Invent.
  Math. \textbf{81} (1985), 449--457.

\bibitem[Ge]{Ge}
A.Gentimis.
Asymptotic dimension of finitely presented group, preprint 2005.

\bibitem[Gr87]{Gr87}
M.~Gromov.
\newblock Hyperbolic groups.
\newblock In {\em Essays in group theory}, 75--263. Springer, New York,
  1987.

\bibitem[Gr93]{Gr93}
M.~Gromov, \emph{Asymptotic invariants of infinite groups},
Geometric Group
  Theory, London Math. Soc. Lecture Note Ser. (G.~Niblo and M.~Roller, eds.),
  no. 182, 1993.

\bibitem[Ha]{Ha}
W.~J. Harvey.
\newblock Boundary structure of the modular group.
\newblock In {\em Riemann surfaces and related topics: Proceedings of the 1978
  Stony Brook Conference},
  245--251, 1981. Princeton Univ. Press.

\bibitem[JS]{JS}
Tadeusz Januszkiewicz, Jacek Swiatkowski.
Filling invariants in systolic complexes and groups.
preprint 2005.



\bibitem[Man]{Man}
J.Manning, 
Fox Geometry of pseudocharacters.  Geom. Topol.  9  (2005), 1147--1185.

\bibitem[MaMi]{MaMi}
H.Masur, Y.N.Minsky. 
Geometry of the complex of curves. I. Hyperbolicity.  Invent. Math.  138  (1999),  no. 1, 103--149.

\bibitem[Mu]{Mu}
J.R.Munkres, "Topology", 2nd Ed. Prentice Hall 2000,

\bibitem[Sch]{Sch}
S.Schleimer.
The end of the curve complex.
preprint, 2006.
math.GT/0608505

\bibitem[Sta68]{Sta}
J.~Stallings, \emph{On torsion free groups with infinitely many ends}, Ann. of
  Math. \textbf{88} (1968), 312--334.


\end{thebibliography}
\end{document}